\RequirePackage{ifpdf}
\ifpdf 
\documentclass[pdftex]{sigma}
\else
\documentclass{sigma}
\fi

\newcommand{ \Ga}{ \Gamma}
\newcommand{ \La}{ \Lambda}
\newcommand{ \la}{ \lambda}

\def\Hom{\mathop{\rm Hom}}
\def\ad{\mathop{\rm ad}}

\begin{document}

\allowdisplaybreaks

\renewcommand{\thefootnote}{$\star$}

\renewcommand{\PaperNumber}{059}

\FirstPageHeading

\ShortArticleName{Projections of Singular Vectors of Verma Modules}

\ArticleName{Projections of Singular Vectors of Verma Modules\\ over Rank 2 Kac--Moody Lie Algebras\footnote{This paper is a
contribution to the Special Issue on Kac--Moody Algebras and Applications. The
full collection is available at
\href{http://www.emis.de/journals/SIGMA/Kac-Moody_algebras.html}{http://www.emis.de/journals/SIGMA/Kac-Moody{\_}algebras.html}}}

\Author{Dmitry FUCHS and Constance WILMARTH}

\AuthorNameForHeading{D. Fuchs and C. Wilmarth}

\Address{Department of Mathematics, University of California, One Shields Ave., Davis CA 95616, USA}
\Email{\href{mailto:fuchs@math.ucdavis.edu}{fuchs@math.ucdavis.edu}, \href{mailto:wilmarth@math.ucdavis.edu}{wilmarth@math.ucdavis.edu}}

\ArticleDates{Received June 29, 2008, in f\/inal form August 24,
2008; Published online August 27, 2008}

\Abstract{We prove an explicit formula for a projection of singular vectors in the Verma module over a rank~2 Kac--Moody Lie algebra onto the universal enveloping algebra of the Heisenberg Lie algebra and of $sl_{2}$ (Theorem~\ref{theorem3}). The formula is derived from a more general but less explicit formula due to Feigin, Fuchs and Malikov
[{\it Funct. Anal. Appl.} {\bf 20} (1986), no.~2, 103--113].
In the simpler case of  $ \mathcal{A}_{1}^{1}$ the formula was obtained in
[Fuchs D., {\it Funct. Anal. Appl.} {\bf 23} (1989), no.~2, 154--156].}


\Classification{17B67}

\section{Introduction}\label{section1}

Let $ \mathcal{G}(A)$ be the complex Kac--Moody Lie algebra corresponding to an $n \times n$ symmetrizable Cartan
matrix $A$, let $N_{-}, H, N_{ + } \subset\mathcal{G}(A)$ be subalgebras generated by the groups of standard generators: $f_{ i }$, $h_i$, $ e_{ i }$, $ i = 1, \ldots, n$.  Then $ \mathcal{G}(A) = N_{ - } \oplus H \oplus N_{ + }$ (as a vector space), the Lie algebras $N_{ - }$ and $N_{ + }$ are virtually nilpotent, and $H$ is commutative. Let $ \la: H \longrightarrow \mathbb{C}$ be a linear functional and let $M$ be a $ \mathcal{G}(A)$-module.  A non-zero vector $w \in M$ is called a singular vector of type $ \la$ if $gw = 0$ for $g \in N_{ + }$ and $hw = \la (h) w$ for $h \in H.$  Let
\begin{gather*} J_{ \la } = \{ \alpha \in \mathcal{U}(N_{ - } ) | \, \exists  \textrm{ a } \mathcal{G}(A)\textrm{-module }M \textrm{ and a singular vector } w \in M\\
\phantom{J_{ \la } =\{}{}\textrm{ of type } \la \textrm{ such that } \alpha w = 0 \}.
\end{gather*}
Obviously $J_{ \la }$ is a left ideal of $ \mathcal{U}(N_{ - })$.  It has a description in terms of Verma modules $M( \la ).$

Let $I_{ \la }$ be a one-dimensional $ (H \oplus N_{ + })$-module with $hu = \la (h)u$, $gu = 0$ for $g \in N_{ + }$ and arbitrary $u \in I_{ \la }.$  The Verma module $M( \la )$ is def\/ined as the $ \mathcal{G}(A)$-module induced by $I_{ \la }$; as a~$ \mathcal{U}(N_{ - })$-module, $M( \la )$ is a free module with one generator $u$; this ``vacuum vector'' $u$ is, with respect to the $ \mathcal{G}(A)$-module structure, a singular vector of type $ \la.$  It is easy to see that $M( \la )$ has a unique maximal proper submodule and this submodule $L( \la )$ is, actually, $J_{ \la }u.$

This observation demonstrates the fundamental importance of the following two problems.
\begin{enumerate}\itemsep=0pt
\item
For which $ \la $ is the module $M( \la )$ reducible, that is, $L( \la ) \ne 0 $?
\item
If $M( \la )$ is reducible, then what are generators of $L( \la )$ (equivalently, what are generators of $J_{ \la }$)?
\end{enumerate}

Problem~1 is solved, in a very exhaustive way, by Kac and Kazhdan~\cite{kk:struct}.
They describe a~subset $ \mathcal{S} \subseteq  H^ {* }$ such that the module $M( \la )$ is reducible if and only if $ \la \in \mathcal{S}$; this subset is a~countable union of hyperplanes.  (See a precise statement in Section~\ref{section2} below.) Actually, $\lambda\in\mathcal{S}$ if and only if $M(\lambda)$ contains a singular vector not proportional to $u$.

A formula for such a singular vector in a wide variety of cases  is given in the work of Feigin, Fuchs and Malikov~\cite{Feigin}.
This formula is short and simple, but it involves the generators $f_{ i }$ raised to complex exponents; when reduced to the classical basis of $ \mathcal{U}(N_{ - } )$ the formula becomes very complicated (as shown in~\cite{Feigin}
 in the example $ \mathcal{G}(A) = sl_{n}).$  There remains a hope that the projection of these singular vectors onto reasonable quotients of $ \mathcal{U}(N_{ - })$ will unveil formulas that possess a more intelligible algebraic meaning, and this was shown to be the case by Fuchs with the projection over
the algebra $ \mathcal{A}_{1}^{1}$ into
$ \mathcal{U}(sl_{2})$ and $\mathcal{U}(\mathcal{H}),$ where $\mathcal H$ is the Heisenberg algebra~\cite{Fuchs}, work which took its inspiration from the earlier investigation of Verma modules over the Virasoro algebra by Feigin and Fuchs~\cite{FeigF}.

In this note we extend these results by providing projections to
$\mathcal{U}(sl_{2})$ and $\mathcal{U}(\mathcal{H})$ of the singular vectors over the family of Kac--Moody Lie algebras $ \mathcal{G}(A)$ of rank 2 (see Theorem~\ref{theorem3} in Section~\ref{section4} and a discussion in Section~\ref{section5}). As in \cite{Fuchs} and \cite{FeigF}, our formulas express the result in the form of an explicit product of polynomials of degree~2 in ${\mathcal U}({\mathcal H})$ and ${\mathcal U}(sl_2)$.

  It is unlikely that this work can be extended to algebras of larger rank.

\section{Preliminaries}\label{section2}
Let $A = (a_{ij} )$ be an integral $n \times n$ matrix with $a_{ij}=2$ for $i=j$ and $a_{ij}\le0$ for $i\ne j$.
We assume that that $ A $ is symmetrizable, that is, $DA = A^{\rm sym},$ where $D = [d_1, \ldots , d_{n} ]$ is diagonal, $d_{i} \ne 0,$ and $ A^{\rm sym} $ is symmetric. To $A$ is associated a  \emph{Kac--Moody} Lie algebra $\mathcal{G}(A)$ def\/ined in the following way.

$\mathcal{G}(A)$ is a complex Lie algebra with the generators $e_i$, $h_i$, $f_i$, $i=1,\dots,n$ and the relations $[h_i,h_j]=0$, $[h_i,e_j]=a_{ij}e_j$, $[h_i,f_j]=-a_{ij}f_j$, $[e_i,e_j]=\delta_{ij}h_i$, $(\ad e_i)^{-a_{ij}+1}e_j=0$, $(\ad f_i)^{-a_{ij}+1}f_j=0$. There is a vector space direct sum decomposition $\mathcal{G}(A)=N_-\oplus H\oplus N_+$ where $N_-,H,N_+\subset\mathcal{G}(A)$ are subalgebras generated separately by $\{f_i\}$, $\{h_i\}$, $\{e_i\}$. Actually, $H$~is a commutative Lie algebra with the basis $\{h_i\}$. We introduce in $H$ a (possibly, degenerate) inner product by the formula $\langle h_i,h_j\rangle=d_ia_{ij}$.

Fix an auxiliary $n$-dimensional complex vector space $T$ with a basis $\alpha_1,\dots,\alpha_n$; Let $\Gamma$ denote a lattice generated by  $\alpha_1,\dots,\alpha_n$, and let $\Gamma_+$ be the intersection of $\Gamma$ with the (closed) positive octant. For an integral linear combination $\alpha=\sum_{i=1}^nm_i\alpha_i$, denote by $G_\alpha$ the subspace of $\mathcal{G}(A)$ spanned by monomials in $e_i$, $h_i$, $f_i$ such that for every $i$, the dif\/ference between the number of occurrences of $e_i$ and $f_i$ equals $m_i$. If $\alpha\ne 0$ and $G_\alpha\ne0$, then $\alpha$ is called a {\em root} of $\mathcal{G}(A)$.  Every root is a positive, or a negative, integral linear combination of $\alpha_i$; accordingly the root is called positive or negative (and we write $\alpha>0$ or $\alpha<0$). Obviously, $N_+=\oplus_{\alpha>0}G_\alpha$, $N_-=\oplus_{\alpha<0}G_\alpha$. Remark that Verma modules have a natural grading by the semigroup $\Gamma_+$.

 For $ \alpha = \sum k_{i} \alpha_{i},$ let $h_\alpha = \sum k_{i} d_{i}^{ -1} h_{i}.$  We can carry the inner product from $H$ to $T$ using the formula $\langle\alpha,\beta\rangle=\langle h_\alpha,h_\beta\rangle$. If $\langle\alpha,\alpha\rangle\ne0$, then we def\/ine a ref\/lection $s_\alpha\colon H^\ast\to H^\ast$ by the formula
 \[
 (s_\alpha\lambda)(h)=\lambda(h)-\frac{2\lambda(h_\alpha)}{\langle\alpha,\alpha\rangle}\langle h_\alpha,h\rangle.
 \]
 The similar formula
 \[
 s_\alpha\beta=\beta-\frac{2\langle\alpha,\beta\rangle}{\langle\alpha,\alpha\rangle}\alpha
 \]
 def\/ines a ref\/lection $s_\alpha\colon T\to T$. ``Elementary ref\/lections'' $s_i=s_{\alpha_i}$ generate the action of  the {\em Weyl group} $W(A)$ of $\mathcal{G}(A)$ in $H^\ast$ and in $H$. In $H^\ast$ we consider, besides the ref\/lections $s_\alpha$ the ref\/lections $s_\alpha^\rho$, $s_\alpha^\rho(\la)=s_\alpha(\lambda+\rho)-\rho$ where $\rho\in H^\ast$ is def\/ined by the formula $\rho(h_i)=1$, $1\le i\le n$.

 The Kac--Kazhdan criterion for reducibility of Verma modules $M( \la )$, mentioned above, has precise statement:

 \begin{theorem}\label{theorem1}
$M( \la )$ is reducible if and only if for some positive root $\alpha$ and some positive integer~$m$,
\begin{equation}
( \la + \rho )(h_{ \alpha }) -\frac{m}{2} \langle\alpha, \alpha\rangle=0.
\label{eq1}
\end{equation}
 Moreover, if $\lambda$ satisfies this equation for a unique pair $\alpha$, $m$, then all non-trivial singular vectors of $M(\lambda)$ are contained in $M(\lambda)_{m\alpha}$.
\end{theorem}

For $m$ and $ \alpha$ satisfying this criterion, Feigin, Fuchs and Malikov \cite{Feigin} give a description for the singular vector  of degree $m\alpha$ in $M( \la ).$ In the case when $\alpha$ is a real root, that is, $\langle\alpha,\alpha\rangle\ne0$, their description is as follows. Let $s_\alpha=s_{i_N}\cdots s_{i_1}$ be a presentation of $s_\alpha\in W(A)$ as a product of elementary ref\/lections. For $\lambda\in H^\ast$, set $\lambda_0=\lambda$, $\lambda_j=s_{i_j}(\la_{j-1}+\rho)-\rho$ for $0<j\le N$. Obviously, the vector $\overrightarrow{\la_{j-1}\la_j}$ is collinear to $\alpha_{i_j}$ (or, rather, to $\langle\alpha_{i_j},\ \rangle$); let $\overrightarrow{\la_{j-1}\la_j}=\gamma_j$. Let $\alpha$ satisf\/ies (for some $m$) the equation~\eqref{eq1}. Then
\begin{equation*}
F(s_\alpha;\lambda)u\qquad \mbox{where}\quad F(s_\alpha;\la)=f_{i_N}^{-\gamma_N}\cdots f_{i_1}^{-\gamma_1} 
\end{equation*}
is a singular vector in $M(\lambda)_{m\alpha}.$ Notice that the exponents in the last formula are, in general, complex numbers. It is explained in \cite{Feigin} why the expression for $F(s_\alpha;\la)$  still makes sense.

\section{The case of rank two}\label{section3}

 In the case $n=2$, a (symmetrizable) Cartan matrix given by
 \[ A =\left( \begin{array}{cc} 2 & {-q} \\   {-p} & 2 \end{array} \right),
 \]
 where $p>0$, $q>0$. Since for $pq\le3$, the algebra $\mathcal{G}(A)$ is f\/inite-dimensional, we consider below the case when $pq\ge4$.

Simple calculations show that $s_{\alpha_{1}}(\alpha_1) = -\alpha_1$, $s_{\alpha_{1}}(\alpha_2) = p\alpha_1 + \alpha_2$, $s_{\alpha_2}(\alpha_1) = \alpha_1+ q \alpha_2$, and $ s_{\alpha_2}(\alpha_2) = -\alpha_2$, and it is easy to check that the orbit of the root $(1,0)$ lies in the curve $qx^2 -pqxy + py^2= q$ and  the orbit of $(0,1)$ lies in $qx^2 - pqxy + py^2 = p.$ (If $pq<4$, these two curves are hyperbolas sharing asymptotes, in the (degenerate) case of $pq=4$, they are pairs of parallel lines with a slope of $\frac q2$.)

  Def\/ine a sequence recursively by $a_0 = 0$, $a_1 = 1,$ and $a_{n} = sa_{n-1} - a_{n-2}$ where $s^2 = pq.$  Then for $\sigma^2 =  \frac{q}{p}$  we can calculate
 \begin{gather*}
 (1,0)   =  (a_1, \sigma a_0),   \\
 s_2((1,0))   =   (1,q) = (a_1, \sigma a_2),   \\
 s_1s_2((1,0))  = (pq-1,q)= (a_3, \sigma a_2),   \\
 s_2 s_1 s_2((1,0))  =  (pq-1, q(pq-2)) = (a_3, \sigma a_4).
 \end{gather*}
 More generally, the following is true.
\begin{proposition}\label{proposition1}
 The orbit of $(1,0)$ consists of points
 \[ \cdots  (a_{2n-1}, \sigma a_{2n-2}), (a_{2n-1}, \sigma a_{2n}), (a_{2n+1}, \sigma a_{2n}), (a_{2n+1}, \sigma a_{2n+2}) \cdots \]
determined by the sequence $ \{a_{n} \}$ above; while the orbit of $(0,1)$ consists of points
 \[  \cdots(\sigma^{-1} a_{2n-2}, a_{2n-1}),
  (\sigma^{-1} a_{2n}, a_{2n-1}),
 (\sigma^{-1} a_{2n}, a_{2n+1}),
 (\sigma^{-1} a_{2n+2}, a_{2n+1}) \cdots \]
 for $n \ge 1.$
 \end{proposition}
 \begin{proof}
 The proof is by induction on $n$.
 \end{proof}

 Obtaining explicit coordinates for the real orbits is straightforward in the af\/f\/ine case, because of the simpler geometry.  For $pq > 4$ an explicit description of the sequence $ \{ a_{n} \}$ is possible using an argument familiar to Fibonnaci enthusiasts:

  \begin{proposition}\label{proposition2}
  The  $nth$ term is
  \begin{gather*}
 a_{n}  =  \frac {1}{  \sqrt {pq - 4}}  \left(  \frac{ \sqrt {pq} + \sqrt {pq - 4}}{2}  \,\right) ^{n} - \frac{1}{ \sqrt{pq - 4}} \left( \frac{ \sqrt{pq} - \sqrt{pq - 4}}{2} \,\right) ^{n}.
   \end{gather*}
  \end{proposition}
  \begin{proof} Direct computation. \end{proof}

 With these real roots now labelled by the sequence $ \{ a_{n} \} ,$ we present the singular vectors indexed by them in the Verma modules over $\mathcal{G}(A)$.  Write $ \lambda = x \lambda_1 + y \lambda_2$ where $ \lambda_{i} (h_{j}) = \delta_{ij},$ so that $ \lambda (h_{1}) = x$ and $ \lambda (h_{2}) = y.$  Let us def\/ine the numbers $\Ga^k_1$, $\Ga^k_2$ by the formulas:  \begin{gather*} \Ga^{2m}_1 = {q\sum_{i=0}^{m-1} (-1)^{i} {2m-i-1 \choose 2m-2i-1} (pq)^{m-i-1}}, \\
 \Ga^{2m}_2= {\sum_{i=0}^{m-1} (-1)^{i} {2(m-1)-i \choose 2(m-1)-2i} (pq)^{m-i-1} }, \\
 \Ga^{2m+1}_1 = {\sum_{i=0}^{m} (-1)^{i} {2m-i \choose 2m-2i} (pq)^{m-i}},\\
 \Ga^{2m+1}_2 ={p\sum_{i=0}^{m-1} (-1)^{i} {2m-i-1 \choose 2m-2i-1} (pq)^{m-i-1}}. \end{gather*}
(Note that  $\Ga^0_1=\Ga^0_2= 0.)$

The formula from \cite{Feigin}  takes in our case the following form.

\begin{theorem}\label{theorem2}
For  the algebra $  \mathcal{G} ( A)$ with Cartan matrix
\begin{gather*}
A =
\left( \begin{array}{cc}
2 & -p \\
-q & 2
\end{array} \right)
\end{gather*}
the singular vectors are as follows:

1. For the root $\alpha=(a_{2n-1},\sigma a_{2n-2}),$ with $m \in \mathbb{N},$ and $ t \in \mathbb{C}$ arbitrary,
\begin{gather*} F(s_{ \alpha}; \lambda)=
f_1^{ \frac{ \Ga_1^{4n - 3}m\strut}{a_{2n-1}} + \Ga_2^{2n - 1}t}f_2^{ \frac{ \Ga_1^{4n - 4}m\strut}{a_{2n-1}} + \Ga_2^{2n - 2}t} f_1^{  \frac{ \Ga_1^{4n - 5}m\strut}{a_{2n-1}} + \Ga_2^{2n - 3}t}\\
\phantom{F(s_{ \alpha}; \lambda)=}{} \cdots f_2^{ \frac{ \Ga_1^{2n}m\strut}{a_{2n-1}} + \Ga_2^2t} f_1^{m} f_2^{ \frac{ \Ga_1^{2n-2}m\strut}{a_{2n-1}} - \Ga_2^2t} \cdots  f_2^{ \frac{ \Ga_1^{2}m\strut}{a_{2n-1}} - \Ga_2^{2n-2}t} f_1^{ \frac{ \Ga_1^{1}m\strut}{a_{2n-1}} - \Ga_2^{2n-1}t}
\end{gather*}
and the vector $F(s_\alpha,\lambda)u$ is singular in $M(\la)= {M\left( \frac{m - \Ga_1^{ 2n - 1 }}{ \Ga_1^{2n-1}} - \Ga_2^{2n-1}t , \Ga_1^{2n-1}   t - 1\right)}$.

2. For $ \alpha = (a_{2n-1}, \sigma a_{2n})$,{\samepage
 \begin{gather*} F(s_{ \alpha} ; \lambda) =
f_2^{ \frac{ \Ga_2^{4n}m\strut}{\sigma^{-1} a_{2n}} + \Ga_2^{2n}t}f_1^{ \frac{ \Ga_2^{4n-1}m\strut}{\sigma^{-1}a_{2n}} + \Ga_2^{2n-1}t}  f_2^{ \frac{ \Ga_2^{4n-2}m\strut}{\sigma ^{-1}a_{2n}} + \Ga_2^{2n-2}t}\\
\phantom{F(s_{ \alpha} ; \lambda) =}{} \cdots f_2^{ \frac{ \Ga_2^{2n + 2}m\strut}{\sigma^{ -1} a_{2n}} + \Ga_2^{2}t}f_1^{ m} f_2^{ \frac{ \Ga_2^{2n}m\strut}{\sigma^{-1}a_{2n}} -  \Ga_2^{2}t} \cdots  f_1^{ \frac{ \Ga_2^3m\strut}{\sigma^{-1}a_{2n}} - \Ga_2^{2n-1}t}f_2^{ \frac { \Ga_2^2m\strut}{\sigma^{-1}a_{2n}} - \Ga_2^{2n}t} \end{gather*}
 and the vector $F(s_\alpha,\lambda)u$ is singular in $M(\la)= {M \left( \Ga_2^{2n} t - 1, \frac{m - \Ga_2^{2n}}{ \Ga_2^{2n}} - \Ga_1^{2n} t\right)}$.}

3. For $ \alpha = (\sigma^{-1}a_{2n-2}, a_{2n-1})$,
\begin{gather*}  F( s_{ \alpha}, \lambda)=
f_2^{ \frac{ \Ga_2^{4n-2}m\strut}{a_{2n-1}} + \Ga_1^{2n-2}t}f_1^{ \frac{ \Ga_2^{4n-3}m\strut}{a_{2n-1}} + \Ga_1^{2n-3}t}  f_2^{ \frac{ \Ga_2^{4n-4}m\strut}{a_{2n-1}} + \Ga_1^{2n-4}t}\\
\phantom{F( s_{ \alpha}, \lambda)= }{} \cdots f_1^{ \frac{ \Ga_2^{ 2n + 1 }m\strut}{a_{2n-1}} + \Ga_1^{1}t}f_2^{ m} f_1^{  \frac{ \Ga_2^{ 2n - 1}m\strut}{a_{2n-1}} - \Ga_1^{1}t} \cdots f_1^{ \frac{ \Ga_2^3m\strut}{a_{2n-1}} - \Ga_1^{2n-3}t}f_2^{ \frac{ \Ga_2^2m\strut}{a_{2n-1}} - \Ga_1^{2n-2}t}
\end{gather*}
 and the vector $F(s_\alpha,\lambda)u$ is singular in $M(\la)= {M\left( \Ga_2^{2n-1} t - 1, \frac{m - \Ga_2^{ 2n - 1 }}{ \Ga_2^{2n-1}} - \Ga_1^{2n-1} t\right)}.$

4. For $ \alpha = (\sigma^{-1}a_{2n}, a_{2n-1})$,
\begin{gather*}  F( s_{ \alpha} ; \lambda) =
f_1^{ \frac{ \Ga_1^{4n-1}m\strut}{\sigma a_{2n}} + \Ga_1^{2n-1}t}f_2^{ \frac{ \Ga_1^{4n-2}m\strut}{\sigma a_{2n}} + \Ga_1^{2n - 2}t} f_1^{ \frac{ \Ga_1^{4n-3}m\strut}{\sigma a_{2n}}  + \Ga_1^{2n-3}t}\\
\phantom{F( s_{ \alpha} ; \lambda) =}{} \cdots f_1^{ \frac{ \Ga_1^{2n + 1}m\strut}{\sigma a_{2n}} + \Ga_1^{1}t}f_2^{m}f_1^{ \frac{ \Ga_1^{ 2n - 1}m\strut}{\sigma a_{2n}}  - \Ga_1^{2n-3}t} \cdots f_2^{ \frac{ \Ga_1^2m\strut}{\sigma a_{2n}} - \Ga_1^{2n-2}t} f_1^{ \frac{ \Ga_1^1m\strut}{\sigma a_{2n}} - \Ga_1^{2n-1}t}
\end{gather*} and the vector $F(s_\alpha,\lambda)u$ is singular in $M( \lambda) = \left( \frac{m - \Ga_1^{2n}}{ \Ga_1^{2n}} - \Ga_2^{2n} t, \Ga_1^{2n} t - 1\right).$
\end{theorem}
\begin{proof}
It must be checked that the vectors given above actually correspond to the Feigin--Fuchs--Malikov (FFM) procedure for obtaining singular vectors, and also that the  Kac--Kazhdan criterion for reducibility is satisf\/ied.  For $ \lambda = x \lambda_1 + y \lambda_2$ and the ref\/lection $s_{ \alpha} = s_{ i_{N}} \cdots s_{ i_{1}} $ (pro\-duct of simple ref\/lections) the algorithm requires successive application of the transformations $s_{1}^{ \rho} := s_1 (\lambda + \rho) - \rho$ and $s_2^{ \rho} := s_2 (\lambda + \rho) - \rho.$
One generates the list
\begin{gather*}
\la^0  = x \la_1 + y \la_2, \qquad
\la^{j}  = s_{ i_{ j } } ( \la^{ j - 1 } + \rho ) - \rho
\end{gather*}
and the auxiliary sequence $ \{ \overrightarrow{ \la^{ j - 1 } \la^{ j }} \}_{ j \ge 1}. $
The algorithm then gives
\[ F(s_{ \alpha }; \la ) = f_{ i_{N}}^{ \theta_{N}} \cdots f_{ i_{ 1 }}^{ \theta_{ 1 }}, \]
where $ \overrightarrow{ \la^{ j - 1 } \la^{ j }} = - \theta_{ j } \alpha_{ i_{ j }}$ (here $\alpha_{ i_{ j }}$
is the functional $\langle h_{ \alpha_{ i_{ j }}}, \cdot\rangle).$

So we f\/irst need to know the decomposition of $s_{ \alpha } $ into elementary ref\/lections for $ \alpha$ in the orbit of $(1,0)$ or $(0,1).$
 Let $S_{i} (m)$ denote the word in $H^{*}$ beginning and ending with $s_{i}$, and containing $m$ $s_{i}$'s.  For example, $S_{1}(3) = s_{1} s_{2} s_{1} s_{2} s_{1}.$
  \begin{lemma}\label{lemma1}
  For real $ \alpha$ as above, $s_{\alpha}$ is the word
  \begin{gather*}
  (a_{2n-1}, \sigma a_{2n-2})  \longleftrightarrow S_1(2n-1), \\
(a_{2n-1}, \sigma a_{2n})  \longleftrightarrow S_2 (2n), \\
(\sigma^{-1}a_{2n}, a_{2n-1})  \longleftrightarrow S_1(2n), \\
(\sigma^{-1}a_{2n-2}, a_{2n-1})  \longleftrightarrow S_2 (2n-1).
 \end{gather*}
 \end{lemma}
 \begin{proof}
 This is an easy induction on $n.$
\end{proof}

The coef\/f\/icients of collinearity $ \theta_{ j }$ have the following description.
\begin{lemma}\label{lemma2}
For $ \la^{k} = \La_1^{k} \la_1 + \La_2^{k} \la_2$ we have
\begin{gather*}
(i) \ \overrightarrow{ \la^{2n+1} \la^{2n+2}}= -( \La_2^{2n+1}+1)\langle h_{ \alpha_2}, \cdot\rangle,\\
(ii) \ \overrightarrow{ \la^{2n} \la^{2n+1}} = -( \La_1^{2n}+1)\langle h_{ \alpha_1}, \cdot\rangle.
\end{gather*}
\end{lemma}
\begin{proof}
One easily computes that $\langle h_{ \alpha_1}, \cdot\rangle = 2 \lambda_1 + q \lambda_2$, $\langle h_{ \alpha_2}, \cdot\rangle = -p \lambda_1 + 2 \lambda_2,$ and further that
\begin{gather*}
s_1^{ \rho} (x \la_1 + y \la_2) = (-x-2) \la_1 + (y + q(x+1)) \la_2,
\\
s_2^{ \rho} (x \la_1 + y \la_2) = (x + p(y+1)) \la_1 + (-y-2) \la_2.
\end{gather*}
Then for $(i)$ it is verif\/ied that $\overrightarrow{ \la^1 \la^2} = -(y+qx+q+1)(-p \la_1 + 2 \la_2) = -( \La_2^1 +1)\langle h_{ \alpha_2}, \cdot\rangle.$  For $n>0,$
\begin{gather*}
 \la^{2n+2}   = s_2^{ \rho} (s_1^{ \rho} ( \la^{2n}))  =s_2^{ \rho}( (- \La_1^{2n}-2) \la_1 + ( \La_2^{2n} + q \La_1^{2n}+q) \la_2) \\
 \phantom{\la^{2n+2}}{}  =(- \La_1^{2n}-2+p \La_2^{2n}+pq \La_1^{2n}+pq+p) \la_1+ (- \La_2^{2n} - q \La_1^{2n} - q -2) \la_2
 \end{gather*}
 while
 \begin{gather*}
 \la^{2n+1} = s_1^{ \rho}( \la^{2n})=(- \La_1^{2n}-2) \la_1+ ( \La_2^{2n}+q \La_1^{2n} +q) \la_2.
 \end{gather*}
 So $ \overrightarrow{ \la^{2n+1} \la^{2n+2}}= -( \La_2^{2n} + q \La_1^{2n} +q+1)( -p \la_1+2 \la_2).
 $   Since $ \la^{2n+1}= s_1^{ \rho}( \la^{2n})=( - \La_1^{2n}-2) \la_1+( \La_2^{2n}+q
 \La_1^{2n} +q) \la_2,$ we have $ - \La_2^{2n+2}-1= -( \La_2^{2n}+q \La_1^{2n}+q+1)$ as desired.
 The argument for~$(ii)$ is similar.
 \end{proof}

 Let us put $\Ga^k=\Ga^k_1x+\Ga^k_2y$. We will also need
 \begin{lemma}\label{lemma3}
\begin{gather*}
(i) \ \Ga^{2n+1} = p \Ga^{2n}- \Ga^{2n-1},\\
(ii) \ \Ga^{2n+2}=q \Ga^{2n+1}- \Ga^{2n}.
\end{gather*}
\end{lemma}
\begin{proof}
These can be verif\/ied directly.
\end{proof}

 We are now in a position to show by induction that the FFM-exponents correspond to the~$ \Ga^{k}$ in the statement of Theorem~\ref{theorem2}.  It suf\/f\/ices by the second lemma to show that
 \begin{gather*}
 \Ga^{2n+1} = \La_1^{2n} + 1
 \qquad \mbox{and}\qquad
 \Ga^{2n+2} = \La_2^{2n+1} + 1.
 \end{gather*}
 Making a change of variable $x+1 \rightarrow x$ and $y+1 \rightarrow y$ one can calculate that
 \begin{gather*}
 \overrightarrow{\la^{0} \la^{1}}  = -x\langle h_{\alpha_1}, \cdot\rangle, \\
 \overrightarrow{ \la^{1} \la^{2}}  = -(qx+y)\langle h_{ \alpha_2}, \cdot\rangle, \\
 \overrightarrow{ \la^{2} \la^{3}}  = -((pq-1)x + py)\langle h_{ \alpha_{1}}, \cdot\rangle, \\
\dots\dots\dots \dots\dots\dots\dots\dots\dots
 \end{gather*}

The FFM exponents are just the coef\/f\/icients of these functionals with reversed sign.  For the base case of our induction, $n=0,$ observe that $ \la^{0} = (x-1) \la_1 + (y-1) \la_2,$ hence $ \La_1^{0} + 1=x= \Ga^{1}$ (using the binomial def\/inition of $ \Ga^{1}),$ while $ \La_2^{1} +1 =y+qx$ since $ s_1^{ \rho} ( \la^{0} ) = (-x-1) \la_1 + (y-1+qx) \la_2,$ agreeing with the binomial sum $ \Ga^{2} = y+qx.$

 Inductively assume that for some $(n-1) > 0$
 \begin{gather*}
 \Ga^{2(n-1)+1} = \La_1^{2n-2} + 1
 \qquad \mbox{and}\qquad
 \Ga^{2(n-1)+2} = \La_2^{2n-1} +1.
 \end{gather*}
 We need to show that
 \begin{gather*}
 \Ga^{2n+1}= \La_1^{2n} +1
 \qquad \mbox{and}\qquad
 \Ga^{2n+2} =  \La_2^{2n+1} +1.
 \end{gather*}
 Just observe that
 \begin{gather*}
 \La_1^{2n} +1  = \La_1^{2n-1} + p( \La_2^{2n-1} +1) +1
  = ( \La_1^{2n-1} +1) + p( \La_2^{2n-1} +1) \\
 \phantom{\La_1^{2n} +1 }{} = (- \La_1^{2n-2} -2 +1)  + p \left[ \La_2^{2n-2} + q( \La_1^{2n-2} +1) +1 \right] \\
\phantom{\La_1^{2n} +1 }{}  = -( \La_1^{2n-2} +1) +p \left[ ( \La_2^{2n-2} +1) +q( \La_1^{2n-2} +1) \right] \\
\phantom{\La_1^{2n} +1 }{} = - \Ga^{2n-1} + p( \La_2^{2n-1} +1)
 = - \Ga^{2n-1} +p \Ga^{2n}
  = \Ga^{2n+1},
 \end{gather*}
 where the last equality comes from Lemma~\ref{lemma3} and the inductive hypothesis is used in the preceding two lines.  The same tack proves that $ \Ga^{2n+2}= \La_2^{2n+1} +1.$

 We now know that the FFM-exponents are as given Theorem~\ref{theorem2}.  It only remains to check that the Kac--Kazhdan criterion (Theorem~\ref{theorem1}) is satisf\/ied. For $m$ and $ \alpha$ satisfying this criterion, \cite{Feigin}~give the prescription for the singular vector  $F(s_{ \alpha }; \la )u$ of degree $m \alpha$ in $M( \la )$; so we need to verify the existence of such $ \alpha$ and $m.$

  For $ \alpha = a \alpha_1 + b \alpha_2$ and $ \la = x \la_1 + y \la_2$, $h_{ \alpha}= ad_1^{-1} h_1 + bd_2^{-1}= \frac{a}{p} h_1 + \frac{b}{q} h_2$, so the criterion can be restated as
 $2(x \la_1 + y  \la_2 + \rho) {\left( \frac{a}{p}h_1 + \frac{b}{q}h_2\right)}   =m\langle a \alpha_1 + b \alpha_2, a \alpha_1+b \alpha_2\rangle. $
 After the calculations this is \[2\left(x \frac{a}{p} +y \frac{b}{q} + \frac{a}{p} + \frac{b}{q}\right) =m\left(\frac{2a^2}{p}-2ab+\frac{2b^2}{q}\right)\] or \[(x+1) \frac{a}{p}+(y+1) \frac{b}{q} = m\left( \frac{a^2}{p} -ab+ \frac{b^2}{q}\right).\]
 So after change of variable $ x+1 \rightarrow x$ and $ y+1 \rightarrow y$ the Kac--Kazhdan criterion becomes \[x \frac{a}{p} + y \frac{b}{q} = m\left( \frac{a^2}{p} - ab+ \frac{b^2}{q}\right).\]

 We show that the integral exponent of the centermost element in the singular vectors in the statement of the theorem precisely meets the integrality requirement of the criterion.  This is a case by case check, and somewhat tedious and technical; let us verify it for roots of type $(a_{2n+1}, \sigma a_{2n}),$ whose singular vector comprises $2n+1$  $f_1$'s and $2n$ $f_2$'s raised to appropriate powers; the centermost exponent is then the $2n+1$-st coef\/f\/icient of collinearity in the FFM procedure, or what we have called $ \Ga^{2n + 1}.$

 A remark, a lemma, and a corollary will show that $\Ga^{2n + 1}$ does what it is supposed to.

\begin{remark}  $ q\Ga_2^{2n+1} = p\Ga_1^{2n},$ as is transparent from the def\/initions of
 $ \Ga^{2n}$ and $ \Ga^{2n+1}.$
 \end{remark}

 The next lemma will relate the root sequence $ \{ a_{n} \}$ to the exponents of the singular vectors.
 \begin{lemma}\label{lemma4}
 The following is true for $ n \ge 0$:
 \begin{gather*}
 \Ga_1^{2n+1} = a_{2n+1},  \qquad
 \Ga_1^{2n+2}  = \sigma a_{2n+2}.
 \end{gather*}
 \end{lemma}
 \begin{proof}
 Induction on $n.$
 \end{proof}

 \begin{corollary}\label{corollary1}
 $ \Ga^{2n+1} = a_{2n+1}x + \frac{p}{q} \sigma a_{2n}y.$
 \end{corollary}
 \begin{proof}
 $ \Ga^{2n+1} = \Ga_1^{2n+1}x + \Ga_2^{2n+1}y
 = a_{2n+1}x + \Ga_2^{2n+1}y
 = a_{2n+1}x +  \frac{p}{q} \Ga_1^{2n} y
 = a_{2n+1}x +  \frac{p}{q}ua_{2n}y.$
 \end{proof}

 Finally the Kac--Kazhdan criterion for $( a_{2n+1}, ua_{2n})$ is
 \begin{gather*}
 a_{2n+1} \frac{x}{p} + ua_{2n} \frac{y}{q} = m\left(\frac{(a_{2n+1})^2}{p} - ua_{2n}a_{2n+1} + \frac{ (ua_{2n})^2}{q}\right)
 \end{gather*}
 or equivalently
 \begin{gather*}
 a_{2n+1}x+ \frac{p}{q}ua_{2n}y   = m((a_{2n+1})^2 -pua_{2n}a_{2n+1}+ \frac{p}{q}(ua_{2n})^2) \\
 \phantom{a_{2n+1}x+ \frac{p}{q}ua_{2n}y}{}  = m( (a_{2n+1})^2 - pua_{2n}a_{2n+1} + \frac{p}{q} \frac{q}{p}(a_{2n})^2) = m' \in \mathbb{N}
 \end{gather*}
 since $ a_{2n+1}$ and $ua_{2n}$ are integral (polynomial in $p$ and $q$).  But the left-hand side here is by the corollary exactly the exponent of the centermost letter in the singular vector, which by the formula
 given in the theorem is an integer; so the Kac--Kazhdan criterion is indeed satisf\/ied in this case.  One can check in similar fashion that the remaining three cases also f\/it the integrality requirement.

The singular vectors in the statement of the theorem appear as follows.  The exponent of the centermost vector in all four cases must be integral: setting this expression in $x$ and $y$ equal to~$m$ one then solves for $x$ (or $y$) in terms of $m$ and $y$ (respectively, $x$); $t$ is then introduced as a~scalar multiple of $y$ (resp., $x$) to minimize notational clutter.
 This completes the proof of the theorem.
 \end{proof}

\section{Projections}\label{section4}

 We next obtain projections of the singular vectors into the Heisenberg algebra, where they factor as products.  While the theorem gives a simple and perhaps the most natural expression for the singular vectors in terms of the $ \Ga^{k}$ a  change of variable is advantageous in the projection and factoring of these vectors in the Heisenberg algebra.  In each case this involves setting the exponent of the vector immediately to the left of the centermost letter equal to a complex variab\-le~$ \alpha$ (which will then depend on $n$).   For example the root $\gamma  = (u^{-1}a_{4}, a_{3}) = (p(pq-2), pq-1)$ has from the theorem the corresponding singular vector:
\begin{gather*} f_1^{ \frac{((pq)^3-5(pq)^2+6pq-1)m}{q(pq-2)}+(pq-1)t} f_2^{ \frac{q((pq)^2-4pq+3)m}{q(pq-2)}+qt} f_1^{ \frac{(pq)^2-3pq+1)m}{q(pq-2)}+t}f_2^{m}\\
\qquad{}\times f_1^{ \frac{(pq-1)m}{q(pq-2)}-t}f_2^{ \frac{qm}{q(pq-2)}-qt}f_1^{ \frac{m}{q(pq-2)}-(pq-1)t}.
\end{gather*}
Taking
\begin{gather*}
\alpha = \frac{( (pq)^2-3pq+1)m}{q(pq-2)} + t
\end{gather*}
the singular vector becomes
\begin{gather*}
f_1^{(pq-1) \alpha -pm}f_2^{q \alpha -m}f_1^{ \alpha}f_2^{m}f_1^{pm- \alpha}f_2^{q(pm - \alpha) -m}f_1^{(pq-1)(pm- \alpha)-pm}
\end{gather*}
which can in turn be rewritten in terms of the $ \Ga^{k} $ as:
\begin{gather*}
f_1^{ \Ga_2^{4} \alpha - \Ga_2^{3}m}f_2^{ \Ga_1^{2} \alpha - \Ga_1^{1}m}f_1^{ \Ga_2^{2} \alpha -  \Ga_2^{1}m}f_2^{m}f_1^{- \Ga_2^{2} \alpha + \Ga_2^{3}m}f_2^{- \Ga_1^{2} \alpha + \Ga_1^{3}m}f_1^{ - \Ga_2^{4} \alpha + \Ga_2^{5}m}.
\end{gather*}
\begin{proposition}\label{proposition3}
Under this change of variable the singular vectors take the following form:

1. For $\alpha=(a_{2n-1}, \sigma a_{2n-2})$ the corresponding $F( s_\alpha; \la)$  is
\begin{gather*}
f_1^{ \Ga_2^{ 2n  - 1} \alpha - \Ga_2^{ 2n - 2 }m } f_2^{ \Ga_1^{ 2n - 3} \alpha - \Ga_1^{ 2n - 4}m}f_1^{ \Ga_2^{2n - 3} \alpha - \Ga_2^{2n - 4}m}f_2^{ \Ga_1^{2n - 5} \alpha - \Ga_1^{2n - 6}m}f_1^{ \Ga_2^{2n - 5} \alpha - \Ga_2^{2n - 6}m}\\
\qquad{}\cdots f_2^{ \Ga_1^{3} \alpha - \Ga_1^{2} m}f_1^{ \Ga_2^{3} \alpha - \Ga_2^{2}m}f_2^{ \Ga_1^{1} \alpha - \Ga_1^{0}m}f_1^{m}f_2^{ - \Ga_1^{1} \alpha + \Ga_1^{2}m}f_1^{- \Ga_2^{3} \alpha + \Ga_2^{4}m}f_2^{ -\Ga_1^{3} \alpha + \Ga_1^{4}m} \\
\qquad{}\cdots  f_1^{ - \Ga_2^{2n-3} \alpha + \Ga_2^{2n-2}m}f_2^{ - \Ga_1^{2n-3} \alpha + \Ga_1^{2n-2}m}f_1^{ - \Ga_2^{2n-1} \alpha + \Ga_2^{2n}m}.
\end{gather*}

2. For $ \gamma = ( a_{2n-1}, ua_{2n})$ the corresponding singular vector $F(s_{\gamma}; \la)$  is
\begin{gather*}
f_2^{ \Ga_1^{2n-1} \alpha - \Ga_1^{2n-2}m} \big( F(s_{(a_{2n-1}, ua_{2n-2})} ; \la \big) f_2^{- \Ga_1^{2n-1} \alpha + \Ga_1^{2n}m}.
\end{gather*}

3. For $ \gamma = (u^{-1}a_{2n-2}, a_{2n-1})$ the singular vector $ F(s_{ \gamma}; \la)$ becomes
\begin{gather*}
 f_2^{ \Ga_1^{2n-2} \alpha - \Ga_1^{2n-3}m}f_1^{ \Ga_2^{2n-2} \alpha - \Ga_1^{2n-3}m}f_2^{ \Ga_1^{2n-4} \alpha - \Ga_1^{2n-5}m}f_1^{ \Ga_2^{2n-4} \alpha - \Ga_2^{2n-5}m} \cdots \\
\qquad{} f_2^{ \Ga_1^{2} \alpha - \Ga_1^{1}m}f_1^{ \Ga_2^{2} \alpha - \Ga_2^{1}m}f_2^{m}f_1^{ - \Ga_2^{2} \alpha + \Ga_2^{3}m}f_2^{ - \Ga_1^{2} \alpha + \Ga_1^{3}m}  \cdots f_1^{ -\Ga_2^{2n-2} \alpha + \Ga_2^{2n-1}m}f_2^{- \Ga_1^{2n-2} \alpha + \Ga_1^{2n-1}m}.
\end{gather*}

4. For $ \gamma = (u^{-1}a_{2n}, a_{2n-1})$ the singular vector is
\begin{gather*}
f_1^{ \Ga_2^{2n} \alpha - \Ga_2^{2n-1}m} \big( F(s_{(u^{-1}a_{2n-2}, a_{2n-1)}} ; \la)\big) f_1^{- \Ga_2^{2n} \alpha + \Ga_2^{2n+1}m}.
\end{gather*}
\end{proposition}

\begin{proof}
This can be established by induction on the number of pairs transformed.
\end{proof}

We now project the singular vectors into the universal enveloping algebra of the three-dimensional Heisenberg algebra $\mathcal H$.  Recall that this is generated by $f_1$, $f_2,$ with $[f_1, f_2] =: h$, $[f_1, h] = [f_2, h] = 0.$ Thus, the projection ${\mathcal U}(N_-)\to{\mathcal U}({\mathcal H})$ is the factorization over the (two-sided) ideal generated by $[f_1,h]$, $[f_2,h]$.

Let $H_{u} := f_{2}f_{1} + uh$ for $u \in  \mathbb{C}.$  Observe that $ H_{u}  H_{v} = H_{v}  H_{u}$, $u,v \in \mathbb{C}.$  The following relations also hold in the Heisenberg (for positive integers, and hence for arbitrary complex numbers $\alpha$, $\beta$, $u$).
\begin{lemma}For $ \alpha, \beta, u \in \mathbb{C}$,
\begin{enumerate}\itemsep=0pt
\item[1)]
$f_2^{ \alpha}  H_{u} = H_{u - \alpha} f_2^{ \alpha}$;
\item[2)]
$f_1^{ \beta} H_{u} = H_{u+ \beta} f_1^{ \beta}$;
\item[3)]
$f_1^{ \alpha} f_2^{n} f_1^{n - \alpha} =  H_{ \alpha }  H_{  \alpha - 1 } \cdots H_{  \alpha - (n-1) }$;
\item[4)]
$f_2^{ \alpha} f_1^{n} f_2^{n - \alpha} = H_{1 - \alpha} H_{2 - \alpha} \cdots H_{n - \alpha}$.
\end{enumerate}
\end{lemma}

\begin{proof}
The calculations follow readily from the complex binomial formula given in \cite{Feigin}: for $ g_{1}, g_{2} \in \mathcal{G}$ a Lie algebra, $ \gamma_{1}, \gamma_{2} \in \mathbb{C},$ we have
\begin{gather*}
g_1^{ \gamma_1} g_2^{ \gamma_2} = g_2^{ \gamma_2} g_1^{ \gamma_1} + \sum_{j_1 = 1}^{ \infty} \sum_{j_2 = 1}^{ \infty} \binom{ \gamma_1}{j_1} \binom{ \gamma_2}{j_2} Q_{j_1 j_2}( g_1, g_2) g_2^{ \gamma_2 - j_2} g_1^{ \gamma_1 - j_1},
\end{gather*}
where $ {\binom{ \gamma}{j} =  \frac{ \gamma^{(j)} }{j!}}$ with $ \gamma^{(j)} = \gamma ( \gamma - 1) \cdots  (\gamma - j + 1)$ and the Lie polynomials $Q_{j_1 j_2}$ can be calculated explicitly, for example using the recursion
\begin{gather*}
Q_{j_1 j_2}(g_1,g_2) = [g_1, Q_{j_1 - 1, j_2}(g_1, g_2)] + \sum_{v = 0}^{ j_2 - 1} \binom{j_2}{v} Q_{j_1 - 1,v}(g_1,g_2) [ \underbrace{ g_2, \dots , [g_2, [g_2 }_{j_2 - v}, g_1]] \dots ]
\end{gather*}
with $Q_{00} = 1$ and $Q_{0,v} = 0$ for $v > 0.$
\end{proof}

Now set, for $r,s \in \mathbb{N}$
\begin{gather*}
 \mathcal{H}_0^{r,s} = 1 \qquad
\textrm{(the empty product)}, \\
 \mathcal{H}_1^{r,s} = \prod_{k=1}^{( \Ga_2^{2} - \Ga_1^{0})m}H_{( - \Ga_2^{2}+ \Ga_2^{3}- \cdots \pm\Ga_2^{r})\alpha+ ( -\Ga_1^{1}+ \Ga_1^{2}- \cdots \pm \Ga_1^{s})m \, +k}.
\end{gather*}
 The for $j \ge 1$ def\/ine
\begin{gather*}
 \mathcal{H}_{2j}^{r,s} =  \prod_{k=1}^{( \Ga_1^{2j} - \Ga_2^{2j})m} H_{ ( \Ga_2^{2j+1}- \Ga_2^{2j+2}+ \cdots \pm \Ga_2^{r}) \alpha + ( -\Ga_1^{2j-1} + \Ga_1^{2j}- \cdots \pm \Ga_1^{s})m - \, (k-1)}, \\
 \mathcal{H}_{2j+1}^{r,s} = \prod_{k=1}^{( \Ga_2^{2j+2}- \Ga_1^{2j})m}H_{(- \Ga_2^{2j+2}+ \Ga_2^{2j+3}- \cdots \pm \Ga_2^{r}) \alpha + ( \Ga_1^{2j} - \Ga_1^{2j+1}+ \cdots \pm \Ga_1^{s})m \, +k}.
\end{gather*}
We will also need the following, not dissimilar, but warranting its own notation:
\begin{gather*}
 \mathcal{ \tilde{H}}_0^{r,s} = 1, \\
 \mathcal{ \tilde{H}}_1^{r,s} = \prod_{k=1}^{( \Ga_2^{2} - \Ga_2^{1})m} \mathcal{H}_{( \Ga_1^{1}- \Ga_1^{2} + \cdots \pm \Ga_1^{r}) \alpha + ( \Ga_2^{2} - \Ga_2^{3} + \cdots \pm \Ga_2^{s})m - (k-1)}.
\end{gather*}

For $j \ge 1$ def\/ine
\begin{gather*}
 \mathcal{ \tilde{H}}_{2j}^{r,s} = \prod_{k=1}^{( \Ga_2^{ 2j + 1 } - \Ga_2^{2j})m} H_{( - \Ga_1^{2j} + \Ga_1^{2j+1} - \cdots \pm \Ga_1^{r}) \alpha + ( \Ga_2^{2j} - \Ga_2^{2j+1} + \cdots \pm \Ga_2^{s})m \, +k}, \\
 \mathcal{ \tilde{H}}_{ 2j + 1 }^{r,s} =  \prod_{k=1}^{ ( \Ga_2^{ 2j + 2 } - \Ga_2^{ 2j + 1 }) m} H_{ ( \Ga_1^{2j+1} - \Ga_1^{2j+2} + \cdots \pm \Ga_1^{r}) \alpha + ( - \Ga_2^{2j+1} + \Ga_2^{2j+2} - \cdots \pm \Ga_2^{s})m - (k-1)}.
\end{gather*}

\begin{theorem}\label{theorem3}
The singular vectors whose words $F(s_{ \alpha}; \la)$ were given in the preceding theorem project to the Heisenberg algebra as the following products:
\begin{enumerate}\itemsep=0pt
\item
The singular vector corresponding to $(a_{2n-1}, ua_{2n-2})$ projects to
\begin{gather*}
\prod_{w=1}^{2n-2} \mathcal{H}_{w}^{2n-1, 2n-3}\, f_1^{( \Ga_1^{2n-1}- \Ga_1^{2n-2})m}.
\end{gather*}
\item
The singular vector corresponding to $(a_{2n-1}, ua_{2n})$ projects to
\begin{gather*}
\prod_{w=1}^{2n-1} \mathcal{H}_{w}^{2n, 2n-2} \, f_2^{( \Ga_1^{2n} - \Ga_1^{2n-1})m}.
\end{gather*}
\item
The singular vector for $(u^{-1}a_{2n-2}, a_{2n-1})$ projects to
\begin{gather*}
\prod_{w=1}^{2n-2} \mathcal{ \tilde{H}}_{w}^{2n-2, 2n-2} f_2^{( \Ga_2^{2n} - \Ga_2^{2n-1})m}.
\end{gather*}
\item
The singular vector for $( u^{-1}a_{2n}, a_{2n-1})$ projects to
\begin{gather*}
\prod_{w=1}^{2n-1} \mathcal{ \tilde{H}}_{w}^{2n-1, 2n-1} f_1^{( \Ga_2^{2n+1} - \Ga_2^{2n})m}.
\end{gather*}
\end{enumerate}
\end{theorem}
\begin{proof}
Induction on $n.$  This is completely straightforward using the change of variables for the singular vectors given in the proposition above.
\end{proof}

\section{Other projections }\label{section5}
One might ask whether projections into other algebras, for instance into $ \mathcal{U}(sl_{2}),$ are equally possible.
A homomorphism from the universal enveloping algebra of $N_{-}$ to $ \mathcal{U}(sl_2)$ is def\/ined when $p>1$ and  $q>1$. It is the factorization over the (two-sided) ideal generated by $[h,f_1]-f_1$ and $[h,f_2]+f_2$ (where, as before, $h=[f_1,f_2]$). Set, for $ u \in \mathbb{C},$
\begin{gather*}
J_{u} = f_2f_1 + uh - \frac{u(u - 1)}{2}.
\end{gather*}
It can then be checked that for $ \beta \in \mathbb{C}$ we have the following in $ \mathcal{U}(sl_2)$ (putting $e = f_1$, $f = f_2):$
\begin{gather*}
f_1^{ \beta} J_{u}  = J_{u + \beta} f_1^ { \beta},\qquad
f_2^{ \beta}J_{u}  = J_{ u - \beta}f_2^{ \beta}
\end{gather*}
as well as, for $n \in \mathbb{N},$
\begin{gather*}
f_1^{ \beta} f_2^{n}f_1^{ n - \beta}  = J_{ \beta} J_{ \beta - 1} \cdots J_{ \beta - (n - 1)}, \\
f_2^{ \beta} f_1^{n} f_2^{ n - \beta}  = J_{ 1 - \beta} J_{ 2 - \beta} \cdots J_{ n - \beta}, \\
f_1^{n} f_2^{ n}  = J_{1} \cdots J_{n}.
\end{gather*}
These properties permit the formal manipulations that af\/ford the factorization results we have already detailed.  Simply substitute $J$ for $H$ in the statement of the projection theorem.

\section{Concluding remarks }\label{section6}
In
\cite{Fuchs}
it is observed that information about singular vectors of Verma modules can be used to obtain information about the homologies of nilpotent Lie algebras.  Namely, the dif\/ferentials of the Bernstein--Gel'fand--Gel'fand resolution $BGG$ of $\mathbb C$ over $N_-$ (see~\cite{Bern}) are presented by matrices whose entries are singular vectors in the Verma modules. Thus, if $V$ is an $N_{ - }$-module that is trivial over the kernel of the projection of $ \mathcal{U}(N^{ - })$ onto $ \mathcal{U}(H)$ or $\mathcal{U}(sl_{2})$ considered above, then our formulas give an explicit description of $BGG \otimes V$ and
$\Hom(BGG,V)$.

\pdfbookmark[1]{References}{ref}
\LastPageEnding

\end{document}